\documentclass[12pt,oneside,reqno]{amsart}
\usepackage{amsmath,amssymb,amsfonts,latexsym,txfonts,pxfonts,wasysym,amsthm}
\usepackage{amssymb,longtable}
\usepackage{hyperref}
\usepackage{lscape,graphicx}%
\usepackage{amsmath}%
\usepackage{amsfonts}%
\usepackage{amssymb}%
\usepackage{graphicx}
\usepackage{times}
\usepackage{cite}

\allowdisplaybreaks \newmuskip\pFqmuskip
\newcommand*\pFq[6][8]{%
  \begingroup
  \pFqmuskip=#1mu\relax
    \mathcode`\.=\string"8000
   \begingroup\lccode`\~=`\.
  \lowercase{\endgroup\let~}\pFqdot
   {}_{#2}F_{#3}{\left[\genfrac..{0pt}{}{#4}{#5};#6\right]}%
  \endgroup
  }
\newcommand{\fs}{\frac{1}{2}}
\title[Another Method for Proving a Hypergeometric Generating relation contiguous to that of  Exton]{Another Method for Proving a Hypergeometric Generating relation contiguous to that of Exton}
\author{\textbf{Shantha Kumari, K.$^{1}$,  Prathima, J.$^{2*}$ and  Arjun K. Rathie$^3$}}
\address{$^{1} $ Department of Mathematics, A. J. Institute of Engineering and Technology, Mangalore, Karnataka State, India}
\email{shanthakk99@gmail.com}
\address{$^{2}$ Department of Mathematics, Manipal Institute of Technology, Manipal, Karnataka State, India}
\email {pamrutharaj@yahoo.co.in}
\address{$^3$ Department of Mathematics, School of Mathematical and Physical Sciences, Central University of Kerala,  Tejaswini Hills, Periye P.O., Kasaragod, 671316, Kerala State, India.}
\email {akrathie@cukerala.ac.in} 
\thanks{$^*$ Corresponding Author}
\begin{document}
\begin{abstract}
The aim of this note is to  establish  an interesting  hypergeometric  generating  relation  contiguous to that of Exton by a short method.
\end{abstract}
\maketitle
\section{INTRODUCTION}
\indent 
The generalized hypergeometric functions with $p$ numerator and $q$ denominator parameters is defined by \cite{rainville}
\begin{align} \label{GHS}
_{p}F_{q}\left[\begin{array}[c]{ccccc} \alpha_{1}, & \ldots, & \alpha_{p} &  & \\ 
                                                   &  &  & ; & x  \\ 
																									\beta_{1}, & \ldots, & \beta_{q} &  &  \end{array} \right]   
																			& =\,_{p}F_{q}\left[\alpha_{1},\ldots,\alpha_{p};\beta_{1},\ldots,\beta_{q};x\right]  \nonumber\\
& ={\displaystyle\sum\limits_{n=0}^{\infty}} \frac{(\alpha_{1})_{n}\ldots (\alpha_{p})_{n}}{(\beta_{1})_{n}\ldots(\beta_{q})_{n}}\; \frac{x^{n}}{n!}%
\end{align}
where $(\alpha)_{n}$ denotes the Pochhammer symbol (or the raised factorial or shifted factorial, since $(1)_{n}=n!$) defined  for every $\alpha \in \mathbb{C} $ by
\begin{equation}
(\alpha)_{n}=\left\{\begin{array}{c} \alpha(\alpha+1) \ldots(\alpha+n-1), \quad n\in\mathbb{N}   \\  1, \qquad \qquad \qquad \qquad  \qquad n=0 \end{array}  \right. 
\end{equation}
For a detailed exposition of this function, we refer the standard texts of Rainville\cite{rainville}, Slater\cite{slater} and Exton\cite{exton}.

 With the help of known result\cite[equ.(2)]{bailey1950}; see also \cite[ p.101, equ.(5)]{erdelyihtf1}
\begin{equation}
\pFq{2}{1}{a,\; a+\frac{1}{2}}{\fs}{x^2}= \fs \left[(1+x)^{-2a} + (1-x)^{-2a} \right]
\end{equation}
\newpage
Exton\cite{exton99} in 1999, obtained the following interesting result.
\begin{align} \label{extonsresult}
\sum_{n=0}^\infty & \frac{(d)_n \; \left(d+ \frac{1}{2}\right)_n }{\left(\frac{1}{2}\right)_n \; n! }\; x^{2n} \;\sum_{m=0}^n A_m\; \frac{(-n)_m \;
\left(-n + \frac{1}{2}\right)_m }{m !} \;y^m \nonumber \\
& = \frac{1}{2}(1+x)^{-2d} \; \sum_{n=0}^\infty A_n\; \frac{(d)_n \;\left(d+\frac{1}{2}\right)_n }{n! }\;\left[\frac{x^2y}{(1+x)^2}\right]^n \nonumber \\
& \quad +   \frac{1}{2}(1-x)^{-2d}\; \sum_{n=0}^\infty A_n \; \frac{(d)_n \left(d+\frac{1}{2}\right)_n }{n! }\left[\frac{x^2y}{(1-x)^2}\right]^n
\end{align}
where $A_m$ is the generalized coefficient. Also, as special case, by letting 
\begin{equation*}
A_n = \frac{(a_1)_n \ldots (a_A)_n}{(h_1)_n \ldots (h_H)_n} = \frac{((a))_n}{((h))_n}
\end{equation*}
he deduced the following result:
\begin{align}
\sum_0^\infty & \frac{(d)_n \;\left(d+ \frac{1}{2}\right)_n }{\left(\frac{1}{2}\right)_n \; n! }\; x^{2n} \;
\pFq{A+2}{H}{(a), \;-n, \; -n+\fs}{(h)}{y}  \nonumber \\
&= \frac{1}{2}(1+x)^{-2d}\; \pFq{A+2}{H}{(a), \;d,\; d+\fs \;}{(h)}{\frac{x^2y}{(1+x)^2}} \nonumber \\
&\quad +  \frac{1}{2}(1-x)^{-2d} \;\pFq{A+2}{H}{(a), \; d, \; d+\fs \;}{(h)}{\frac{x^2y}{(1-x)^2}}
\end{align}
In 2000, with the aid of Baileys identity\cite[Equ (3)]{bailey1950}
\begin{equation}\label{baileysidentity}
\pFq{2}{1}{a, \; a +\fs \; }{\frac{3}{2}}{x^2} = \frac{1}{2x (1-2a)}\left[(1+x)^{1-2a}-(1-x)^{1-2a}\right]
\end{equation}
Malani et al.\cite{malani} established the following interesting hypergeometric generating relation by employing the same technique used by Exton. 
\begin{align} \label{malaniresult}
\sum_{n=0}^\infty & \frac{(d)_n \; \left(d+ \frac{1}{2}\right)_n }{\left(\frac{3}{2}\right)_n \; n! }\; x^{2n} \;\sum_{m=0}^n A_m\; \frac{(-n)_m \; \left(-n - \frac{1}{2}\right)_m }{m !} \;y^m \nonumber \\
& = \frac{1}{2x(1-2d)}(1+x)^{1-2d} \; \sum_{n=0}^\infty A_n\; \frac{(d)_n \;\left(d-\frac{1}{2}\right)_n }{n! }\;\left[\frac{x^2y}{(1+x)^2}\right]^n \nonumber \\
& \quad -   \frac{1}{2x(1-2d)}(1-x)^{1-2d}\; \sum_{n=0}^\infty A_n \; \frac{(d)_n \left(d-\frac{1}{2}\right)_n }{n! }\left[\frac{x^2y}{(1-x)^2}\right]^n
\end{align}
Also, as special case, by letting 
\begin{equation*}
A_n = \frac{(a_1)_n \ldots (a_A)_n}{(h_1)_n \ldots (h_H)_n} = \frac{((a))_n}{((h))_n}
\end{equation*}
they deduced the following result :
\newpage
\begin{align}
\sum_{n=0}^\infty & \frac{(d)_n \;\left(d+ \frac{1}{2}\right)_n }{\left(\frac{3}{2}\right)_n \; n! }\; x^{2n} \;
\pFq{A+2}{H}{(a), \;-n, \; -n-\fs}{(h)}{y}  \nonumber \\
&= \frac{1}{2x(1-2d)}\left\{(1+x)^{1-2d}\; \pFq{A+2}{H}{(a), \;d,\; d-\fs \;}{(h)}{\frac{x^2y}{(1+x)^2}} \right. \nonumber \\
&\qquad \qquad -  \left. (1-x)^{1-2d} \;\pFq{A+2}{H}{(a), \; d, \; d-\fs \;}{(h)}{\frac{x^2y}{(1-x)^2}}\right\}
\end{align}
The aim of this short note is to provide another method for proving the  result \eqref{malaniresult} due to Malani, without using the result \eqref{baileysidentity}
\section{Derivation of the result \eqref{malaniresult}}
Inorder to prove the result  \eqref{malaniresult}, we proceed as follows:\\
 
Without loss of generality we can assume that
\begin{align} \label{9}
& \frac{1}{2x}(1+x)^{1-2d} \; \sum_{n=0}^\infty A_n\; \frac{(d)_n \;\left(d-\frac{1}{2}\right)_n }{n! }\;\left[\frac{x^2y}{(1+x)^2}\right]^n \nonumber \\
& \quad -   \frac{1}{2x}(1-x)^{1-2d}\; \sum_{n=0}^\infty A_n \; \frac{(d)_n \left(d-\frac{1}{2}\right)_n }{n! }\left[\frac{x^2y}{(1-x)^2}\right]^n \nonumber \\
& = \sum_{n=0}^\infty a_{2n+1}\; x^{2n}
\end{align}
Then, it is not much difficult to see that the coefficient $a_{2n+1}$ of $x^{2n}$ in the expansion, after some simplification is obtained as 
\begin{equation}\label{a2n+1}
a_{2n+1} = -\sum_{m=0}^n A_m\; \frac{(d)_m \;\left(d- \frac{1}{2}\right)_m }{ m! } \quad \frac{(2d+2m-1)_{2n-2m+1}}{(2n-2m+1)!} \; y^{m} 
\end{equation}
A simple calculation shows that
\begin{equation}
(2d+2m-1)_{2n-2m+1} = \frac{ \Gamma(2d)\; 2^{2n} \; (d)_n \; \left(d+\fs\right)_n }{\Gamma(2d-1)\; 2^{2m} \; (d)_m \; \left(d - \fs \right)_m}
\end{equation}
 and 
\begin{equation}
(2n-2m+1) ! = \frac{2^{2n} \; \left(\frac{3}{2}\right)_n \; n!}{2^{2m} \; (-n)_m \; \left(-n- \fs \right)_m}
\end{equation}
Substituting these values in \eqref{a2n+1}, we get
\begin{equation} a_{2n+1} = (1-2d)\; \frac{(d)_n \;\left(d+ \frac{1}{2}\right)_n }{\left(\frac{3}{2}\right)_n\; n! } \quad \sum_{m=0}^n A_m\;  \frac{(-n)_m  \; \left(-n - \frac{1}{2}\right)_m }{ m! } \; y^{m} 
\end{equation}
Finally substituting this value  in \eqref{9} and dividing both sides of the resulting expression by $(1-2d)$, we get the result \eqref{malaniresult}.
This completes the proof of \eqref{malaniresult}.
\subsection*{Remark :} The result \eqref{malaniresult} due to  Malani et al. was re-derived by Qureshi et al.\cite[Equ. (2.3)]{qureshi2002}, in 2002.

\end{document}